\documentclass[letterpage,12pt]{article}
 \usepackage[english]{babel}
 \usepackage[pdftex]{graphicx}
  \usepackage{wrapfig,multirow,fullpage,floatrow,mathtools}
  \usepackage{fancybox,enumerate, setspace}
  
  \usepackage{sectsty}
  \sectionfont{\large}
  \usepackage{amsfonts,amsmath,amsthm,array,amssymb}

\DeclarePairedDelimiter\floor{\lfloor}{\rfloor}

 \include{latexsymb,color}
 \include{pxfonts}
 \include{displaymath}
 \include{tipa}

\begin{document}

\begin{center}
{\Large \textbf{Fitting UK strike data using a discrete analogue of gamma-Lomax distribution}}
\end{center}

\smallskip

\begin{center}
{\large   Indranil Ghosh$^{1}$,\quad Ayman Alzaatreh$^2$, \quad GG Hamedani$^{3}$ }
\\
\vspace*{0.5cm}
{\large $^{1}$ University of North Carolina, Wilmington, North Carolina, USA \\
$^{2}$ American University of Sharjah, Sharjah, UAE\\
 $^{3}$ Marquette University, Milwaukee, USA.\\
 
 Corresponding author E-mail address: ghoshi@uncw.edu}
\end{center}

\bigskip

\begin{abstract}
This article represents how certain types of blockades in any industrial (heavy industries) production, in particular, industrial strikes can be modeled with the proposed discrete probabilistic distribution as a baseline distribution. We considered the number of outbreaks of strikes in the coal mining industry, the vehicle manufacturing industry, and the transpose industry in the UK obtained from  Consul (1989). We fitted those data sets with the proposed discrete gamma-Lomax distribution and compared the fit with the discrete generalized Pareto distribution (Consul, 1989). For this purpose, we explore the basic properties of the discrete gamma-Lomax distribution including but not limited to: cumulative distribution,
survival, probability mass, quantile and hazard functions, genesis and
rth-order moments; consider maximum likelihood estimation under the normal set up as well as under the censored data set scenario. It is observed that the newly proposed model
 can be useful to describe strikes arising from various types of industries.
\end{abstract}

\bigskip


\section{Introduction}
Since the beginning of the first world war (and may be even prior to that)  there has been a lot of industrial disputes ranging from low labor wages to uncomfortable work conditions that led to
several strikes in the UK. According to the leading daily newspaper in Britain, the Guardian, the number of workers who went on strike in Britain last year fell to the lowest level since the $1890$'s, when Queen Victoria was on the throne. Furthermore, data from the Office for National Statistics in the UK show $33,000$ workers were involved in labour disputes in $2017,$ down from $154,000$ a year earlier. This is the lowest number since records began in $1893,$ the year of Britain's first national coal strikes, when the figure was $634,000.$ The major seminal event that is worthwhile  mentioning here is  the miners' strike of $1984-1985.$ This was a major industrial action to shut down the British coal industry in an attempt to prevent colliery closures. At its height, the strike involved $142,000$ mineworkers. The number of person-days of work lost to the strike was over $26,000,000,$ days making it the largest since the $1926$ general strike (source: Wikipedia).  There are other industries in the UK which have received major set backs in terms of productivity and  loss in revenue because of strikes by it's labor forces due to several reasons.
Needless to say, it has been a matter of great concern to the industry as to how  one can analyze these number of strikes and take appropriate measures. There exists  a sizable number of research  articles where  quantitative insights in strikes are discussed. Among them some noteworthy models are
models from Van der Velden (Velden, $2000$), Skeels and McGrath (Skeels \& McGrath, $1991$),
Leigh (Leigh, $1984$), Buck (Buck, $1984$) and Mauleon and Vannetelbosch (Mauleon \&
Vannetelbosch, $1998$). The first four models calculate strike activity and show the relations between certain variables
and the strike activity. Many variables are used, including future variables, to form the models.
The influence of these variables on the strike activity is calculated by regression analysis. From
this analysis it seems that the  present variables perform better than future variables. Furthermore, it
seems that wage and unemployment are important variables in calculating the strike activity.
The last model, the Mauleon and Vannetelbosch model, uses concepts from game theory to
determine the influence of profit-sharing on the strike activity. The introduction of a profit sharing
scheme increases the strike activity if the bargaining between the union and the
employer takes place at industry level, but reduces the strike activity if the bargaining takes
place at firm level.

\noindent  The main purpose of this article is to establish that the industrial strikes data, specifically strikes data set arising out of the several different industries in the UK, can be described by a discrete probability model, namely the discrete  gamma-Lomax model.  We begin our discussion by providing a general framework which leads to our specific discrete gamma-Lomax distribution from it's continuous model.

Let $F(x)$ be the cumulative distribution function (c.d.f.) of any random variable $X$ and $r(t)$
be the probability density function (p.d.f.) of a random variable $R$ defined on $[0, \infty).$ The c.d.f. of the gamma-$X$ family of distributions defined by Alzaatreh, et al. (2013, 2014) is given by

\begin{equation}
f_{X}(x)
= \frac{1}{\Gamma(\alpha)\beta^{\alpha}}f_{R}(x)\left(-\log\left[1-F_{R}(x)\right]\right)^{\alpha-1}
\left(1-F_{R}(x)\right)^{1/\beta-1}, \quad  x\in \mathbb{R^{+}}.
\end{equation}

\noindent If the random variable $R$ follows the Lomax distribution with the density function $f_{R}(x)=k\theta^{-1}\left(1+\frac{x}{\theta}\right)^{-(k+1)}, \quad x>0$, then (1) reduces to the gamma- Lomax distribution as

\begin{equation}
f_{X}(x)
=\frac{1}{(x+\theta)\Gamma(\alpha)c^{\alpha}}\left\{1+\frac{x}{\theta}\right\}^{-1/c}
\left(\log\left[1+\frac{x}{\theta}\right]\right)^{\alpha-1},
\end{equation}

\noindent $x>0,$  where $c,\alpha,\theta$ are positive parameters.

\noindent Note that if $X$ is replaced by $X +\theta$ , then (2) reduces to the gamma-Pareto distribution which was proposed and studied in Alzaatreh et al. (2012a).

\noindent From (2), the c.d.f. of the gamma-Lomax distribution is given by

\begin{equation}
F_{X}(x)
=\frac{1}{\Gamma(\alpha)}\gamma\left\{\alpha, c^{-1}\log\left[1+\frac{x+1}{\theta}\right]\right\}, \quad x\geq 0,
\end{equation}

\noindent where $\gamma(\alpha,t)=\int_{0}^{t}u^{\alpha-1} e^{-u} du$ is the incomplete gamma function.

\noindent There are a variety of works available in the literature that extends the Lomax distribution under the continuous paradigm. We mention them chronologically as follows.
Note that (1) was studied by Cordeiro et al. (2015) as a particular case of Zografos- Balakrishnan $(G)$ family of distributions, where $G$ is any baseline continuous  distribution (from the perspective of a gamma generated model). Lemonte and Coredeiro. (2013) studied a five-parameter continuous distribution, the so-called McDonald Lomax distribution, that extends the Lomax distribution. Ghitany et al. (2007) studied   properties of a new parametric distribution generated by Marshall and Olkin  extended family of distributions based on the Lomax model.

\noindent However, not much work has been done towards discrete Lomax mixture type models. The most recent reference that the authors can mention here is the work of  Prieto et al. (2014), in which the authors considered the discrete generalized Pareto model (mixing with zero-inflated Poisson distribution) in modeling  road accident blackspots data. It is observed that their discrete model is unimodal and the p.m.f. (probability mass function) is always a decreasing function. So the model is somehow restricted in nature. Now, the proposed discrete gamma-Lomax distribution (henceforth, for short, DGLD) is not  always decreasing. Consequently, it has greater flexibility.

\noindent The discrete gamma-Lomax distribution can be defined in the following way

\begin{eqnarray}
g(x)
&=&P(x\leq dX< x+1)\notag\\
&=&S(x)-S(x+1)\notag\\
&=&\frac{1}{\Gamma(\alpha)}\left[\gamma\left\{\alpha, c^{-1}\log\left[1+\frac{x+1}{\theta}\right]\right\}
-\gamma\left\{\alpha, c^{-1}\log\left[1+\frac{x}{\theta}\right]\right\}\right],
\end{eqnarray}

\noindent  $x\in \mathbb{N^{*}};$  where $\mathbb{N^{*}}=\mathbb{N}\cup \left\{0\right\}$. Also,  $S(x)$  is the survival function corresponding to (4).

\noindent From (4), the c.d.f. and the survival function of the DGLD are given as follows

\begin{equation}
G(x)
=\frac{1}{\Gamma(\alpha)} \gamma\left\{\alpha, c^{-1}\log\left[1+\frac{\floor*{x+1}}{\theta}\right]\right\}, \quad x\geq 0,
\end{equation}

\bigskip

\begin{equation}
S(x)
=1-\frac{1}{\Gamma(\alpha)} \gamma\left\{\alpha, c^{-1}\log\left[1+\frac{\floor*{x+1}}{\theta}\right]\right\}, \quad x\geq 0,
\end{equation}

\noindent where $\floor*{x}=\max \{m\in\mathbb{Z}\mid m\le x\},$ the floor function.
The probability mass function in (4) is useful for maximum likelihood estimation of the parameters, see Section 3. The survival function in (6) is useful for censored maximum likelihood estimation of the parameters, see Subsection 3.2. Figure \ref{fig1} depicts plots for the DGLD where the scale parameter $\theta=1$ and for various values of the shape parameters $c$ and $\alpha.$ The plots indicate that DGLD can have various shapes including reversed J and right skewed unimodal shapes. 


\medskip

\noindent The paper is organized as follows. In Section $2,$ we discuss some structural properties of the proposed DGLD, including shapes, moments, entropy and order statistics. In Section $3,$  maximum likelihood estimation under regular set up as well as under censored data set up are  discussed.  Section 4 deals with certain characterizations of DGLD. Three different real life data sets are considered to illustrate the applicability of the DGLD in Section 5.  Finally, some concluding remarks are provided in Section 6.

\section{Structural properties}
In this section, we discuss some important structural properties of the DGLD. At first, we have the following lemma.

\smallskip

\noindent {\bf Lemma 1.} If a random variable $Y$ follows the gamma-Lomax distribution with parameters $c,\alpha,\theta$,  then $X=\floor*{Y}$ follows the DGLD($c,\alpha,\theta$).

\bigskip
\noindent \textit{Proof.} Follows immediately from (3).

\bigskip

The hazard function associated with DGLD is

\begin{equation}
h(x)
=\frac{\left[\gamma\left\{\alpha, c^{-1}\log\left[1+\frac{x+1}{\theta}\right]\right\}
-\gamma\left\{\alpha, c^{-1}\log\left[1+\frac{x}{\theta}\right]\right\}\right]}{\Gamma(\alpha)-\gamma\left\{\alpha, c^{-1}\log\left[1+\frac{x}{\theta}\right]\right\}}, x\in \mathbb{N^{*}}.
\end{equation}

\bigskip

\noindent Discrete hazard rates arise in several common situations in reliability theory where
clock time is not the best scale on which to describe lifetime. For example, in ammunition
reliability, the number of rounds fired until the first failure (say) is more important than the age in failure.
A similar scenario can be envisioned  when a piece of equipment operates in cycles and the observation is
the number of cycles successfully completed prior to the failure. In other situations a device
is monitored only once per time period and the observation then is the number of time
periods successfully completed prior to the failure of the device (for details, see Shaked et al., 1995).

\bigskip

\noindent Before discussing other properties of the DGLD($c,\alpha,\theta$), we consider the following  series expressions:

\begin{enumerate}

\item
\begin{equation}
\gamma\left\{\alpha, x\right\}=\displaystyle \sum_{m=0}^{\infty}(-1)^{m}\frac{x^{m+\alpha}}{m!(m+\alpha)},
\end{equation}

\noindent  from Nadarajah and Pal (2008).

\item For any $c\in \mathbb{R}$

\begin{equation}
\left\{\log\left[1+\frac{x}{\theta}\right]\right\} ^{c}
=c\sum_{k=0}^{\infty}\binom{k-c}{k}\sum_{j=0}^{\infty} \frac{(-1)^{j}}{(c-j)}\binom{k}{j}P_{k,j}\left(\frac{x}{\theta}\right)^{k+c}
=\phi_{(c,k, j)}\left(\frac{x}{\theta}\right)^{k+c},
\end{equation}
\end{enumerate}

\noindent where  $\phi_{(c,k, j)}=\displaystyle c\sum_{k=0}^{\infty}\binom{k-c}{k}\sum_{k=0}^{\infty} \frac{(-1)^{j}}{(c-j)}\binom{k}{j}$,
$c_{k}=\frac{(-1)^{k+1}}{k+1}$, $P_{j,0}=1$ and $P_{j,k}=\displaystyle\frac{1}{k} \sum_{m=1}^{k}(jm-k+m)c_{m}P_{jk-m},$ $k=1,2,\cdots$.

\bigskip

\noindent Therefore, the p.m.f. in Eq. (4) can be rewritten as
\begin{equation}
P(X=x)=g(x)
=\sum_{m=0}^{\infty}(-1)^{m}\frac{(c^{-1})^{m+\alpha}}{m!(m+\alpha)}
\left\{\phi_{(m+\alpha,k, j)}\left[\left(\frac{x+1}{\theta}\right)^{k+m+\alpha}-\left(\frac{x}{\theta}\right)^{k+m+\alpha}\right]\right\},
\end{equation}

\noindent $x\in \mathbb{N^{*}}.$

\bigskip

\noindent {\bf Theorem 1.} The DGLD has a DFR property for  $\alpha\leq 1$.

\bigskip

\noindent \textit{Proof.}  From (4), it is not difficult to see that the gamma-Lomax distribution possesses a DFR property for $\alpha\leq 1$. Hence the result.

\bigskip

\noindent {\bf Theorem 2.} The DGLD is unimodal and the mode is  $x=m$, where $m\in \left\{ \floor*{x_{0}}-1, \floor*{x_{0}},
\floor*{x_{0}}+1\right\}$ and $x_{0}=\theta\left\{\exp\left(c(\alpha-1)/(c+1)\right)-1\right\}.$ Furthermore, if $\floor{x_{0}}=0,$  then the mode is $m=0$ or 1.

\bigskip

\textit{Proof.} From Theorem $2$ in Alzaatreh et al. (2012a), the gamma-Pareto distribution is unimodal and the mode is at $x_{0} = 0$ or $x_{0} =\theta\exp\left(c(\alpha-1)/(c+1)\right).$ Therefore, the mode of DGLD is at $x_{0} = 0$   or
$x_{0}=\theta\left\{\exp\left(c(\alpha-1)/(c+1)\right)-1\right\}.$ The rest of the proof follows from Theorem $2$ in Alzaatreh et al. (2012b).

\bigskip

\indent Stochastic ordering is an integral tool to judge comparative behaviors  of random variables. Many stochastic orders exist and have various applications. Theorem $3$ and Corollary $1$  (below) give some results on the stochastic orderings of the DGLD. The orders considered here are  the stochastic order $\leq_{st},$  and the expectation order $\leq_{E}.$

\smallskip

\noindent {\bf Theorem 3.} The DGLD($c,\alpha,\theta$) has the following properties.

\begin{itemize}

\item Suppose $X_{1}\sim DGLD(c,\alpha,\theta_{1})$ and  $X_{2}\sim DGLD(c,\alpha,\theta_{2})$. If $\theta_{1}>\theta_{2}$, then $X_{1}\leq_{st}X_{2}$.

\item Suppose $X_{1}\sim DGLD(c_{1},\alpha,\theta)$ and  $X_{2}\sim DGLD(c_{2},\alpha,\theta)$. If $c_{1}>c_{2}$, then $X_{1}\leq_{st}X_{2}$.

\item Suppose $X_{1}\sim DGLD(c,\alpha_{1},\theta)$ and  $X_{2}\sim DGLD(c,\alpha_{2},\theta)$. If $\alpha_{1}>\alpha_{2}$, then $X_{1}\leq_{st}X_{2}$.
\end{itemize}

\bigskip

\textit{Proof.} Follows immediately from the c.d.f. of the DGLD.

\smallskip

 \noindent Next, we describe the expectation ordering in the next Corollary which follows from Theorem $3.$

\newpage

\noindent  {\bf Corollary 1.}
 \begin{itemize}

\item Suppose $X_{1}\sim DGLD(c,\alpha,\theta_{1})$ and  $X_{2}\sim DGLD(c,\alpha,\theta_{2})$. If $\theta_{1}>\theta_{2}$, then $X_{1}\leq_{E}X_{2}$.

\item Suppose $X_{1}\sim DGLD(c_{1},\alpha,\theta)$ and  $X_{2}\sim DGLD(c_{2},\alpha,\theta)$. If $c_{1}>c_{2}$, then $X_{1}\leq_{E}X_{2}$.

\item Suppose $X_{1}\sim DGLD(c,\alpha_{1},\theta)$ and  $X_{2}\sim DGLD(c,\alpha_{2},\theta)$. If $\alpha_{1}>\alpha_{2}$, then $X_{1}\leq_{E}X_{2}$.
\end{itemize}

\subsection{Moments}

\noindent The $r^{\text{th}}$ moment of DGLD is given by

\begin{align}
&E\left(X^{r}\right)\notag\\
&=\frac{1}{\Gamma(\alpha)}\sum_{x=0}^{\infty} x^{r} \left[\gamma\left\{\alpha, c^{-1}\log\left[1+\frac{x+1}{\theta}\right]\right\}
-\gamma\left\{\alpha, c^{-1}\log\left[1+\frac{x}{\theta}\right]\right\}\right]\notag\\
&=\frac{1}{\Gamma(\alpha)}\sum_{x=0}^{\infty} x^{r}\sum_{m=0}^{\infty}(-1)^{m}\frac{(c^{-1})^{m+\alpha}}{m!(m+\alpha)}
\left\{\phi_{(m+\alpha,k, j)}\left[\left(\frac{x+1}{\theta}\right)^{k+m+\alpha}-\left(\frac{x}{\theta}\right)^{k+m+\alpha}\right]\right\}.
\end{align}

\bigskip

\noindent {\bf Theorem 4.} If $c <1/ r,$ then the $r^{\text{th}}$ moment of the DGLD($c,\alpha,\theta$) exists.

\bigskip

\noindent \textit{Proof.} Assume  that $X$ follows DGLD. Now, Alzaatreh et al. (2012a) showed that if $c <1/ r, $ then $E(X^{r})$
exists for all $r$. The rest of the proof follows from the fact that $0 \leq  \floor*{X} \leq X.$

\bigskip

\noindent Let $X$ have a DGLD($c,\alpha,\theta$) distribution. Then, the probability generating function of $X$ can be expressed as

\begin{eqnarray}
G_{X}(s)
&=&E\left(s^{X}\right)\notag\\
&=& \sum_{x=0}^{\infty} s^{x}\sum_{m=0}^{\infty}(-1)^{m}\frac{(c^{-1})^{m+\alpha}}{m!(m+\alpha)}
\left\{\phi_{(m+\alpha,k, j)}\left[\left(\frac{x+1}{\theta}\right)^{k+m+\alpha}-\left(\frac{x}{\theta}\right)^{k+m+\alpha}\right]\right\}.
\end{eqnarray}

\noindent The corresponding moment generating function is

\begin{eqnarray*}
M_{X}(t)
&=&E\left(\exp(tX)\right)\notag\\
&=& \sum_{x=0}^{\infty} \exp(tx)\sum_{m=0}^{\infty}(-1)^{m}\frac{(c^{-1})^{m+\alpha}}{m!(m+\alpha)}
\left\{\phi_{(m+\alpha,k, j)}\left[\left(\frac{x+1}{\theta}\right)^{k+m+\alpha}-\left(\frac{x}{\theta}\right)^{k+m+\alpha}\right]\right\}.
\end{eqnarray*}

\bigskip

\noindent {\bf Theorem 5.} Let $\mu_{[r]}=E\left[X(X-1)\cdots (X-r+1)\right]$ denote the descending $r^{\text{th}}$ order factorial moment.
Then,

$$\mu_{[r]}=\sum_{m=0}^{\infty}(-1)^{m}\frac{(c^{-1})^{m+\alpha}}{m!(m+\alpha)}\left\{\phi_{(m+\alpha,k, j)}\left(r(r-1)\cdots 1\right)
\left[\left(\frac{r+1}{\theta}\right)^{k+m+\alpha}-\left(\frac{r}{\theta}\right)^{k+m+\alpha}\right]\right\} +\mu_{r+1}.$$

\bigskip

\textit{Proof.} Follows immediately by successive differentiation of (10) and then substituting $s=1$.

\subsection{Order statistics}
 Let $X_{1}, X_{2},\cdots, X_{n}$ be a random sample drawn from (4). Then, the probability mass function  of the $i^{th}$ order statistic, $X_{i:n},$ is  given by

\begin{eqnarray*}
P(X_{i:n}=x)
&=&\frac{n!}{(i-1)!(n-i)!}\int_{F(x-1)}^{F(x)} u^{i-1}(1-u)^{n-i}du\notag\\
&=&\frac{n!}{(i-1)!(n-i)!}\sum_{j=0}^{n-i}(-1)^{j}\binom{n-i}{j}\int_{F(x-1)}^{F(x)} u^{i+j-1} du\notag\\
&=&\frac{n!}{(i-1)!(n-i)!}\sum_{j=0}^{n-i}(-1)^{j}\binom{n-i}{j}\frac{1}{(i+j)\left(\Gamma(\alpha)\right)^{i+j}}\notag\\
&&\times \left\{A_{1}(x)-A_{2}(x)\right\},
\end{eqnarray*}

\noindent where

$$A_{1}(x)=\sum_{k_{1}=0}^{\infty}\sum_{k_{2}=0}^{\infty}\cdots \sum_{k_{i+j}=0}^{\infty}
\frac{(-1)^{s_{i+j}\left(c^{-1}\right)^{s_{i+j}+(i+j)\alpha}}}{p_{i+j}}
\left\{\phi_{(\delta_{1},\delta_{2}, s_{i+j}+(i+j)\alpha)}\left(\frac{x}{\theta}\right)^{s_{i+j}+(i+j)\alpha}\right\},$$

and

$$A_{2}(x)=\sum_{k_{1}=0}^{\infty}\sum_{k_{2}=0}^{\infty}\cdots \sum_{k_{i+j}=0}^{\infty}
\frac{(-1)^{s_{i+j}\left(c^{-1}\right)^{s_{i+j}+(i+j)\alpha}}}{p_{i+j}}
\left\{\phi_{(\delta_{1},\delta_{2}, s_{i+j}+(i+j)\alpha)}\left(\frac{x-1}{\theta}\right)^{s_{i+j}+(i+j)\alpha}\right\},$$

\noindent where $s_{i+j}=\sum_{\ell=1}^{i+j}k_{\ell}$, $p_{i+j}=\prod_{\ell=1}^{i+j} k_{\ell}!(k_{\ell}+\alpha)$.
Now, one can use (11) to obtain a general $r^{\text{th}}$ order moment of $X_{i:n}$.

\bigskip

\noindent The distribution of maximum and minimum order statistics, and the distribution of the range can be derived as follows. 

\noindent Let $X_{i}$, $i=1,2,\cdots, n$ be independent DGLD with parameters $c_{i},\alpha_{i}, \delta_{i}$. Define $U=\min (X_{1},X_{2},\cdots, X_{n})$ and $W=\max(X_{1},X_{2},\cdots, X_{n}).$
Then the  c.d.f. of $U,$ from (8), will be

\begin{eqnarray*}
P(U\leq u)
&=&1-\prod_{i=1}^{n}\left\{1-\frac{1}{\Gamma(\alpha_{i})}\sum_{k=0}^{\infty}\frac{(-1)^{k}}{k!(k+\alpha_{i})}
\left[c^{-1}\log\left(\frac{x}{\theta}\right)\right]^{k+\alpha_{i}}\right\}\notag\\
&=&1-\prod_{i=1}^{n}\left\{1-\frac{1}{\Gamma(\alpha_{i})}\sum_{k=0}^{\infty}\frac{(-1)^{k}}{k!(k+\alpha_{i})}
\left\{\phi_{(k+\alpha_{i},m, \alpha)}\left(\frac{u}{\theta}\right)^{m+k+\alpha_{i}}\right\}\right\}.
\end{eqnarray*}

\noindent Hence, the p.m.f. of $U$ is

\begin{align*}
P(U=u)
&=P(U\leq u)-P(U\leq u-1)\notag\\
&=\prod_{i=1}^{n}\left\{1-\frac{1}{\Gamma(\alpha_{i})}\sum_{k=0}^{\infty}\frac{(-1)^{k}}{k!(k+\alpha_{i})}
\left\{\phi_{(k+\alpha_{i},m, \alpha)}\left[\left(\frac{u}{\theta}\right)^{m+k+\alpha_{i}}-\left(\frac{u-1}{\theta}\right)^{m+k+\alpha_{i}}\right]\right\}\right\}.
\end{align*}

\noindent Next, by similar approach, the p.m.f. of $W$ will be

\begin{align*}
P(W=w)
&=P(W\leq w)-P(W\leq w-1)\notag\\
&=\prod_{i=1}^{n}\left\{\frac{1}{\Gamma(\alpha_{i})}\sum_{k=0}^{\infty}\frac{(-1)^{k}}{k!(k+\alpha_{i})}
\left\{\phi_{(k+\alpha_{i},m, \alpha)}\left[\left(\frac{w}{\theta}\right)^{m+k+\alpha_{i}}-\left(\frac{w-1}{\theta}\right)^{m+k+\alpha_{i}}\right]\right\}\right\}.
\end{align*}

\bigskip

\noindent For the distribution of the range ($R$), from Kabe et al. (1969), we may write  $P(X_{1:n}=X_{n:n}=y)=\left\{F(y)-F(y-1)\right\}^{n},$ for $y=0,1,2,\cdots$. If $R=X_{n:n}-X_{1:n}$ denotes the range of the order statistics, then

\begin{equation*}
P(R=0)
= \left\{\frac{1}{\Gamma(\alpha)}\right\}^{n}\left[\sum_{k=0}^{\infty}\frac{(-1)^{k}}{k!(k+\alpha)}
\phi_{(k+\alpha_{i},m, \alpha)}\left\{\left(\frac{y}{\theta}\right)^{m+k+\alpha}-\left(\frac{y-1}{\theta}\right)^{m+k+\alpha}\right\}\right].
\end{equation*}

\subsection{Entropy and Stress-Strength Parameter}
The cumulative residual entropy proposed by Rao et al. $(2004).$ According to them, it is more general
than the Shannon Entropy in that its definition is valid in  both the
continuous and discrete domains, secondly it possesses more general
mathematical properties than the Shannon entropy and thirdly,  it
can be easily computed from the sample data and they converge
asymptotically  to the true values.
It is given by $\eta_{CRE}=-\sum_{x=0}^{\infty}P(X>x)\log P(X>x).$
 In our case  it is given by

 \begin{align*}
\eta_{CRE}
 &= \sum_{j=0}^{\infty}\sum_{x=0}^{\infty}\frac{1}{\left(\Gamma(\alpha)\right)^{j+1}}\notag\\
&\times \left[1-
 \left(\sum_{k=0}^{\infty}\frac{(-1)^{k}}{k!(k+\alpha)}\right)
  \left\{\sum_{k_{1}=0}^{\infty}\sum_{k_{2}=0}^{\infty}\cdots \sum_{k_{j}=0}^{\infty}
\frac{(-1)^{s_{j}\left(c^{-1}\right)^{s_{j}+j\alpha}}}{p_{j}}
\left\{\phi_{(s_{j}+j\alpha,m, \alpha)}\left(\frac{x}{\theta}\right)^{s_{j}+j\alpha+m}\right\}\right\}\right].
\end{align*}

\bigskip

\indent The stress-strength parameter $R = P(X >Y)$ is a measure of component reliability and
its estimation problem, when $X$ and $Y$ are independent and follow a specified common
distribution has been discussed widely in the literature. Suppose that the random variable
$X$ is the strength of a component which is subjected to a random stress $Y.$ Estimation of
$R$ when $X$ and $Y$ are independent and identically distributed following a well-known
distribution, has been considered in the literature. Many applications of the stress strength
model, for its own nature, are related to engineering or military problems. There
are also natural applications in Medicine or Psychology, which involve the comparison
of two random variables, representing for example the effect of a specific drug or
treatment administered to two groups, control and test. Almost all of these studies
consider continuous distributions for $X$ and $Y,$ because many practical applications of
the stress-strength model in engineering fields presuppose continuous quantitative data.
A complete review is available in Kotz et al. (2003). However, in this regard, a relatively
small amount of work is devoted to the discrete or categorical data. Data may be discrete by
nature. For example, the stress pattern in a step-stress accelerated life test can be treated
as a discrete random variable of which the possible values can be obtained from all stress
levels, and the corresponding probabilities can be obtained from the acting times of each
stress level. Moreover, the stress state of a component can be categorized based on the
characteristic of the external loads. The stress-strength parameter, in discrete case, is defined as

$$R=P(X >Y)=\sum_{x=0}^{\infty} f_{X}(x)F_{Y}(x),$$ where $f_{X}(.),$ and $F_{Y}(.)$ denote the p.m.f. and c.d.f. of the independent discrete random variables $X$ and $Y,$ respectively. Next, we have the following result.

\bigskip

\noindent {\bf Theorem 6.} Let $X\sim DGLD(c_{1},\alpha_{1},\theta)$ and $Y\sim DGLD(c_{2},\alpha_{2},\theta)$ be independent random variables. Then, the expression for the  reliability parameter is given by

\begin{eqnarray*}
R&=&\frac{1}{\Gamma(\alpha_{1})\Gamma(\alpha_{2})}\sum_{x=0}^{\infty}\sum_{m_{1}=0}^{\infty}\sum_{m_{2}=0}^{\infty}(-1)^{m_{1}+m_{2}}\frac{(c^{-1}_{1})^{m_{1}+\alpha_{1}}}{m_{1}!(m_{1}+\alpha_{1})}\frac{(c^{-1}_{2})^{m_{2}+\alpha_{2}}}{m_{2}!(m_{2}+\alpha_{2})}\notag\\
&&\times\left\{\phi_{(m_{1},\alpha_{1},k, j)}\phi_{(m_{2},\alpha_{2},k, j)}\left[\left(\frac{x+1}{\theta}\right)^{2k+m_{1}+m_{2}+\alpha_{1}+m_{2}}-
\left(\frac{x}{
\theta}\right)^{k+m_{1}+\alpha_{1}}\left(\frac{x+1}{\theta}\right)^{k+m_{2}}\right]\right\}.
\end{eqnarray*}

\section{Estimation}
\subsection{Regular maximum likelihood estimation}
Several approaches for parameter estimation are proposed in the literature but the maximum likelihood method is the most commonly employed. The maximum likelihood estimators (MLEs) enjoy desirable properties and can be used for constructing confidence intervals and regions and also in test-statistics. The normal approximation for these estimators in large sample distribution theory is easily handled either analytically or numerically.
To apply the method of maximum likelihood for estimating the parameter vector $\Delta=\left(c,\alpha,\theta\right)^{(T)}$ of DGLD distribution, assume that $\vec{x}=(x_1, x_2, \cdots, x_n)^{T}$  is a random sample
of size $n$ from an $X\sim DGLD(c,\alpha,\theta)$ distribution. The log-likelihood function becomes

\begin{equation}
\ell=-n\log \Gamma(\alpha)+\sum_{i=1}^{n}\log\left[\left[\gamma\left\{\alpha, c^{-1}\log\left[1+\frac{x_{i}+1}{\theta}\right]\right\}
-\gamma\left\{\alpha, c^{-1}\log\left[1+\frac{x_{i}}{\theta}\right]\right\}\right]\right].
\end{equation}

\noindent Equation (13) can be maximized using a readily available statistical softwares such as $R$ ({\tt optim} function) or SAS ({\tt PROC NLMIXED}), or by solving the nonlinear likelihood equations obtained by differentiating (13). For interval estimation and hypothesis tests, we can use standard likelihood techniques based on the observed information matrix. For example, the asymptotic covariance matrix of $\Delta$ can be approximated by the inverse observed information matrix evaluated at $\Delta$. Likelihood ratio (LR) tests can be performed for the proposed distribution in the usual way.

\subsection{Censored maximum likelihood estimation}
One may also consider the estimation under censored data set.
Censoring is common in lifetime data sets. There are many types of censoring: type I censoring, type II censoring, and others. A general form known as multi-censoring can be described as follows: there are $m$ lifetimes of which

\begin{itemize}
\item $m_0$ have failed at times $T_1 , \cdots, T_{m_{0}};$
\item $m_1$ have failed at times belonging to $(S_{i -1} , S_{i}]$ ,\quad  $i = 1, \cdots, m_{1};$
\item $m_2$ have survived the times $R_{i}$ , $i = 1 , \cdots,m_{2}$ but no longer observed. It is quite obvious that,
$m = m_0 + m_1 + m_2.$
\end{itemize}

\bigskip

\noindent For the multi-censoring data, the associated log likelihood function will be

\begin{align}
&\log L\left(c,\alpha,\theta\right)\notag\\
&=-\left(m_{0}+m_{1}\right)\log\left(\Gamma(\alpha)\right)
+\sum_{j=1}^{m_{0}}\log\left[\gamma\left\{\alpha, c^{-1}\log\left[1+\frac{T_{j}+1}{\theta}\right]\right\}
-\gamma\left\{\alpha, c^{-1}\log\left[1+\frac{T_{j}}{\theta}\right]\right\}\right]\notag\\
&+\sum_{j=1}^{m_{0}}\log\left[\gamma\left\{\alpha, c^{-1}\log\left[1+\frac{S_{j}+1}{\theta}\right]\right\}
-\gamma\left\{\alpha, c^{-1}\log\left[1+\frac{S_{j}}{\theta}\right]\right\}\right]\notag\\
&+\sum_{j=1}^{m_{0}}\log\left[1-\frac{1}{\Gamma(\alpha)} \gamma\left\{\alpha, c^{-1}\log\left[1+\frac{\floor*{R_{j}+1}}{\theta}\right]\right\}\right].
\end{align}

\bigskip

\noindent The maximum likelihood estimators of $c,\alpha,\theta,$ say $\widehat{c},\widehat{\alpha},\widehat{\theta},$ will be the values of $(c,\alpha,\theta)$ upon maximizing (14).

\section{ Characterizations of the DGLD}

\bigskip

The problem of characterizing a distribution is an important problem in
applied sciences, where an investigator is vitally interested to know if their
model follows the right distribution. \ To this end the investigator relies on
conditions under which their model would follow a specified distribution. In
this section, we present two characterizations of the DGLD based on: $\left(
i\right)$ conditional expectation of certain function of a random variable
and $\left( ii\right)$ the reverse hazard function.

\bigskip

\subsection{ Characterization of the DGLD in terms of the conditional
expectation of certain function of a random variable}

\textbf{Proposition 1.\ } Let \ $X:\Omega\rightarrow%
\mathbb{N}
^{\ast}=%
\mathbb{N}
\cup\left\{  0\right\}  $ \ be a random variable. \ The p.m.f. of \ $X$ \ is
\ $\left(  4\right)  $ \ if and only if

\bigskip%

\begin{align*}
&  E\left\{  \left[  \gamma\left\{  \alpha,c^{-1}\log\left[  1+\frac
{X+1}{\theta}\right]  \right\}  +\gamma\left\{  \alpha,c^{-1}\log\left[
1+\frac{X}{\theta}\right]  \right\}  \right]  \text{ }|\text{ }X\leq k\right\}
\\
&  =\gamma\left\{  \alpha,c^{-1}\log\left[  1+\frac{k+1}{\theta}\right]
\right\}
.\text{\ \ \ \ \ \ \ \ \ \ \ \ \ \ \ \ \ \ \ \ \ \ \ \ \ \ \ \ \ \ \ \ \ \ \ \ \ \ \ \ \ \ \ \ \ \ \ \ \ \ \ \ \ \ \ }%
\left( 15\right)
\end{align*}

\bigskip

\textit{Proof.} \ If \ $X$ \ has p.m.f. \ $\left(  4\right)  $ , then for \ $k\in%
\mathbb{N}
^{\ast},$ the left-hand side of $\left(15 \right)  $ will be

\bigskip%

\begin{align*}
& \left(  G\left(  k\right)  \right)  ^{-1}\sum_{x=0}^{k}\left[  \gamma
^{2}\left\{  \alpha,c^{-1}\log\left[  1+\frac{x+1}{\theta}\right]  \right\}
-\gamma^{2}\left\{  \alpha,c^{-1}\log\left[  1+\frac{x}{\theta}\right]
\right\}  \right] \\
&  =\gamma\left\{  \alpha,c^{-1}\log\left[  1+\frac{k+1}{\theta}\right]
\right\}  .
\end{align*}

\bigskip

Conversely, if \ $\left(  15\right)  $ \ holds, \ then

\bigskip%

\begin{align*}
&\sum_{x=0}^{k}\left\{  \left[  \gamma\left\{  \alpha,c^{-1}\log\left[
1+\frac{x+1}{\theta}\right]  \right\}  +\gamma\left\{  \alpha,c^{-1}%
\log\left[  1+\frac{x}{\theta}\right]  \right\}  \right]  \text{ }g\left(
x\right)  \right\} \\
&=G\left(  k\right)  \gamma\left\{  \alpha,c^{-1}\log\left[  1+\frac
{k+1}{\theta}\right]  \right\}
.\text{\ \ \ \ \ \ \ \ \ \ \ \ \ \ \ \ \ \ \ \ \ \ \ \ \ \ \ \ \ \ \ \ \ \ \ \ \ \ \ \ \ \ \ \ \ \ }%
\left( 16\right)
\end{align*}

\bigskip

From \ $\left( 16\right)  $ , we also have

\bigskip%

\begin{align*}
&  \sum_{x=0}^{k-1}\left\{  \left[  \gamma\left\{  \alpha,c^{-1}\log\left[
1+\frac{x+1}{\theta}\right]  \right\}  +\gamma\left\{  \alpha,c^{-1}%
\log\left[  1+\frac{x}{\theta}\right]  \right\}  \right]  \text{ }g\left(
x\right)  \right\} \\
&  =\left\{  G\left(  k\right)  -g\left(  k\right)  \right\}  \gamma\left\{
\alpha,c^{-1}\log\left[  1+\frac{k}{\theta}\right]  \right\}
.\text{\ \ \ \ \ \ \ \ \ \ \ \ \ \ \ \ \ \ \ \ \ \ \ \ \ \ \ \ \ \ \ \ \ \ \ \ \ }%
\left( 17\right)
\end{align*}

\bigskip

Now, subtracting \ $\left( 17\right)  $ \ from \ $\left(16\right),$
yields

\bigskip%

\begin{align*}
&  \gamma\left\{  \alpha,c^{-1}\log\left[  1+\frac{k+1}{\theta}\right]
\right\}  g\left(  k\right) \\
&  =G\left( k\right)  \left\{  \left[  \gamma\left\{  \alpha,c^{-1}%
\log\left[  1+\frac{k+1}{\theta}\right]  \right\}  -\gamma\left\{
\alpha,c^{-1}\log\left[  1+\frac{k}{\theta}\right]  \right\}  \right]
\right\}  .
\end{align*}

\bigskip

From the above equality, we have

\bigskip%

\[
\frac{g\left(  k\right)  }{G\left(  k\right)  }=1-\left(  \frac{\gamma\left\{
\alpha,c^{-1}\log\left[  1+\frac{k}{\theta}\right]  \right\}  }{\gamma\left\{
\alpha,c^{-1}\log\left[  1+\frac{k+1}{\theta}\right]  \right\}  }\right)  ,
\]

\bigskip

\noindent which is the reverse hazard function of the random variable $X$ with
the p.m.f. \ $\left(4\right).$

\bigskip

\subsection{Characterization of the DGLD in terms of the reverse hazard
function}

\bigskip

\noindent \textbf{Proposition 2.\ }Let \ $X:\Omega\rightarrow%
\mathbb{N}
^{\ast}=%
\mathbb{N}
\cup\left\{  0\right\}  $ \ be a random variable. \ The p.m.f. of \ $X$ \ is
\ $\left(  4\right)  $ \ if and only if its reverse hazard function satisfies
the difference equation

\bigskip%

\begin{align*}
&  r_{G}\left(  k+1\right)  -r_{G}\left(  k\right) \\
=  &  \frac{\gamma\left\{  \alpha,c^{-1}\log\left[  1+\frac{k}{\theta}\right]
\right\}  \gamma\left\{  \alpha,c^{-1}\log\left[  1+\frac{k+2}{\theta}\right]
\right\}  -\gamma^{2}\left\{  \alpha,c^{-1}\log\left[  1+\frac{k+1}{\theta
}\right]  \right\}  }{\gamma\left\{  \alpha,c^{-1}\log\left[  1+\frac
{k+1}{\theta}\right]  \right\}  \gamma\left\{  \alpha,c^{-1}\log\left[
1+\frac{k+2}{\theta}\right]  \right\}  ,\text{
\ \ \ \ \ \ \ \ \ \ \ \ \ \ \ \ \ \ \ \ \ \ \ \ \ \ }\left( 18\right)  }%
\end{align*}

with initial condition $r_{G}\left(  0\right)=1.$

\bigskip

\textit{Proof.} \ If \ $X$ \ has mpf \ $\left( 4\right)$, then clearly \ $\left(
18\right)  $ holds. \ Now, if \ $\left( 18\right)$ holds, then for
\ $x\in%
\mathbb{N}
$, we have

\bigskip%

\begin{align*}
&  \sum_{k=0}^{x-1}\left\{  r_{G}\left(  k+1\right)  -r_{G}\left(  k\right)
\right\} \\
&  =\sum_{k=0}^{x-1}\left\{  \frac{\gamma\left\{  \alpha,c^{-1}\log\left[
1+\frac{k}{\theta}\right]  \right\}  }{\gamma\left\{  \alpha,c^{-1}\log\left[
1+\frac{k+1}{\theta}\right]  \right\}  \text{\ }}-\frac{\gamma\left\{
\alpha,c^{-1}\log\left[  1+\frac{k+1}{\theta}\right]  \right\}  }%
{\gamma\left\{  \alpha,c^{-1}\log\left[  1+\frac{k+2}{\theta}\right]
\right\}  \text{\ }}\right\} \\
&  =-\frac{\gamma\left\{  \alpha,c^{-1}\log\left[  1+\frac{x}{\theta}\right]
\right\}  }{\gamma\left\{  \alpha,c^{-1}\log\left[  1+\frac{x+1}{\theta
}\right]  \right\}  \text{\ }},
\end{align*}

\noindent or%

\[
r_{G}\left(  x\right)  -r_{G}\left(  0\right)  =-\frac{\gamma\left\{
\alpha,c^{-1}\log\left[  1+\frac{x}{\theta}\right]  \right\}  }{\gamma\left\{
\alpha,c^{-1}\log\left[  1+\frac{x+1}{\theta}\right]  \right\}  \text{\ }},
\]

\noindent or, in view of $r_{G}\left(  0\right)  =1$, we have

\bigskip%

\[
r_{G}\left(  x\right)  =1-\left(  \frac{\gamma\left\{  \alpha,c^{-1}%
\log\left[  1+\frac{k}{\theta}\right]  \right\}  }{\gamma\left\{
\alpha,c^{-1}\log\left[  1+\frac{k+1}{\theta}\right]  \right\}  }\right)  .
\]

\section{Application}
In this section, the discrete gamma-Lomax distribution (DGLD) is applied to several data sets. These data sets are taken from Consul (1989). The three data sets represent the observed frequencies of the number of outbreaks of strike in four leading industries in the  U.K. during $1948-1959.$  These industries are coal-mining, vehicle manufacture and Transpose. The data are depicted in Tables $1-3.$ Consul (1989) fitted the data for the four industries to the generalized Poisson distribution (GPD) and the results showed that the GPD does not give an adequate fit to the coal-mining and transpose industries. In this section we show that the discrete gamma-Lomax distribution provide good fit to all data sets.
Tables $1-3$ contain the results of fitting these data sets to the discrete gamma-Lomax and the generalized Poisson distribution. The method of maximum likelihood estimation is used to estimate the model parameters. The chi-square goodness of fit is employed to test the adequacy of the fitted models.

\begin{table}[H]
\centering
\caption{The number of outbreaks of strike in the coal-mining industry in UK}
\begin{tabular}{|c|c|c|c|}
\hline
x-value&	Observed&	DGLD&	GPD\\
\hline
0&	46&	45.99&	50.01\\
\hline
1&	76&	75.59&	65.77\\
\hline
2&	24&	26.23&	32.23\\
\hline
3&	9&	6.32&	7.23\\
\hline
$>=4$	&1&	1.87&	0.76\\
\hline
Total&	156&	156&	156\\
\hline
\end{tabular}
\end{table}

\begin{table}[H]
\centering
\caption{The estimated parameters  and goodness of fit for the outbreaks of strike in the coal-mining industry in UK data. }
\label{app1}
\begin{tabular}{|c|c|c|c|c|}
\hline
Model             & Parameters &$\chi^{2}$& df &p-value   \\ \hline
DGLD       & $\hat{\alpha}=4.5109$  & 1.7279&1&0.1887 \\
            & $\hat{c}=0.0468$& & &\\
            & $\hat{\theta}=6.2260$& & &    \\
            \hline
 GPD             & $\hat{\lambda}=4.5109$  & 4.5234&2&0.0334 \\
            & $\hat{\theta}=1.1377$& & &  \\
 \hline
\end{tabular}
\end{table}

\bigskip

\begin{table}[H]
\centering
\caption{The number of outbreaks of strike in the vehicle-manufacture  industry in UK}
\begin{tabular}{|c|c|c|c|}
\hline
x-value&	Observed&	DGLD&	GPD\\
\hline
0&	110&	109.90&	109.82\\
\hline
1&	33&	33.43&	33.36\\
\hline
2&	9&	8.98&	9.24\\
\hline
3&	3&	2.54&	3.58\\
\hline
$>=4$	&1&	1.15&	\\
\hline
Total&	156&	156&	156\\
\hline
\end{tabular}
\end{table}

\begin{table}[H]
\centering
\caption{The estimated parameters  and goodness of fit for the outbreaks of strike in the  vehicle-manufacture  industry in UK data. }
\label{app1}
\begin{tabular}{|c|c|c|c|c|}
\hline
Model             & Parameters &$\chi^{2}$& df &p-value   \\ \hline
DGLD       & $\hat{\alpha}=4.5109$  & 0.1082&1&0.7422 \\
            & $\hat{c}=0.0468$& & &\\
            & $\hat{\theta}=6.2260$& & &    \\
            \hline
 OGPD             & $\hat{\lambda}=- 0.144$  & 0.06&1&0.8065 \\
            & $\hat{\theta}=0.351$& & &  \\
 \hline
\end{tabular}
\end{table}

\bigskip

\begin{table}[H]
\centering
\caption{The number of outbreaks of strike in the transpose  industry in UK}
\begin{tabular}{|c|c|c|c|}
\hline
x-value&	Observed&	DGLD&	$GPD^{*}$\\
\hline
0&	114&	114.20&	114.41\\
\hline
1&	35&	34.13&	26.01\\
\hline
2&	4&	5.58&	4.83\\
\hline
3&	2&	1.35&	0.85\\
\hline
$>=4$	&1&	0.74& 9.88\\
\hline
Total&	156&	156&	156\\
\hline
\end{tabular}
\end{table}

\begin{table}[H]
\centering
\caption{The estimated parameters  and goodness of fit for the outbreaks of strike in the transpose  industry in UK data. }
\label{app1}
\begin{tabular}{|c|c|c|c|c|}
\hline
Model             & Parameters &$\chi^{2}$& df &p-value   \\ \hline
DGLD       & $\hat{\alpha}=7.6182$  & 0.8707&1&0.3508 \\
            & $\hat{c}=0.1569$& & &\\
            & $\hat{\theta}=0.3166$& & &    \\
            \hline
 $GPD^{*}$             & $\hat{\lambda}=0.098$  & 12.788&2&0.0017 \\
            & $\hat{\theta}=0.31$& & &  \\
 \hline
\end{tabular}
\end{table}

\noindent (*) values in Tables 5 and 6 are different from the computed values in Consul $(1989, p.120)$. The corrected values are computed in this Table.
From Tables $1-4,$ the GPD provides an adequate fit to vehicle-manufacture industry but does not provide an adequate fit to Coal-mining and transpose industries. The DGLD provides adequate fit to all data sets. Furthermore, it can be noted that the DGLD distribution fits the left and right tails of the three data sets well. These data sets show that DGLD can provide an adequate fit to industrial strike data sets. Therefore, the DGLD can be used as a baseline distribution for modeling the number of strikes in industries.

\section{Conclusion}
In this paper we have proposed a new discrete analogue of the continuous gamma-Lomax distribution  and derived some of its distributional properties. The DGLD distribution offers great flexibilities in terms of shapes for the probability mass functions
 and hazard rate functions. The application section shows that DGLD can be useful in fitting various strike data sets related to UK industries.The estimation of the model parameters are discussed in both regular case (with all the data available) and under the multi-censoring set up. From the application section, it appears that DGLD distribution provides a better alternative to the existing GPD probability models. 
 
\section{Declarations} 
\subsection*{Availability of data and material}
All the data sets utilized in this manuscript are available on the specific references made freely. 

\subsection*{Competing interests}
The authors have no conflict of interest.

\subsection*{Funding}
The authors did not receive any funding in preparing this manuscript.

\subsection*{Authors' contributions}
All three authors, Drs. Indranil Ghosh, Ayman Alzaatreh and G.G. Hamedani have equally contributed in this manuscript.

\subsection*{Acknowledgements}
The authors acknowledges several of the references from which some useful ideas have generated in preparing this manuscript.

\end{document}